\pdfoutput=1 
\documentclass[journal]{IEEEtran}

\usepackage{amsmath,graphicx}
\usepackage[ruled]{algorithm2e}
\usepackage{threeparttable}
\usepackage{dcolumn}
\usepackage{booktabs}
\usepackage{amsmath,graphicx}
\usepackage{url}
\usepackage{cite}
\usepackage{psfrag}
\usepackage{amssymb}
\usepackage{amsbsy}
\usepackage{graphicx}
\usepackage{array}
\usepackage{multirow}
\usepackage{bm}
\usepackage{subfigure}
\usepackage{cite}
\usepackage{verbatim}
\usepackage{colortbl}
\usepackage{algorithmic}
\newcommand{\xx}{{\bf x}}
\newcommand{\yy}{{\bf y}}
\newcommand{\ffx}{{\bf f}_x}
\newcommand{\ffy}{{\bf f}_y}
\newcommand{\ccx}{{\bf c}_x}
\newcommand{\ccy}{{\bf c}_y}
\newcommand{\cc}{{\bf c}}
\newcommand{\bbx}{{\bf b}_x}
\newcommand{\bby}{{\bf b}_y}
\newcommand{\bb}{{\bf b}}

\begin{document}

\title{Image Deblurring Using Derivative Compressed Sensing for Optical Imaging Application}

\author{Mohammad~Rostami,~\IEEEmembership{Student~Member,~IEEE,}
        Oleg~Michailovich,~\IEEEmembership{Member,~IEEE,}
        and~Zhou~Wang,~\IEEEmembership{Member,~IEEE}%
\thanks{All the authors are with the Department of Electrical and Computer Engineering, University of Waterloo, N2L 3G1 Ontario.}}

\markboth{IEEE Transactions on Image Processing,~Vol.~XX, No.~XX, January~2012}%
{Shell \MakeLowercase{\textit{et al.}}: Manuscript submission.}

\maketitle

\begin{abstract}
Reconstruction of multidimensional signals from the samples of their partial derivatives is known to be a standard problem in inverse theory. Such and similar problems routinely arise in numerous areas of applied sciences, including optical imaging, laser interferometry, computer vision, remote sensing and control. Though being ill-posed in nature, the above problem can be solved in a unique and stable manner, provided proper regularization and relevant boundary conditions. In this paper, however, a more challenging setup is addressed, in which one has to recover an image of interest from its noisy and blurry version, while the only information available about the imaging system at hand is the amplitude of the generalized pupil function (GPF) along with partial observations of the gradient of GPF's phase. In this case, the phase-related information is collected using a simplified version of the Shack-Hartmann interferometer, followed by recovering the entire phase by means of derivative compressed sensing. Subsequently, the estimated phase can be combined with the amplitude of the GPF to produce an estimate of the point spread function (PSF), whose knowledge is essential for subsequent image deconvolution. In summary, the principal contribution of this work is twofold. First, we demonstrate how to simplify the construction of the Shack-Hartmann interferometer so as to make it less expensive and hence more accessible. Second, it is shown by means of numerical experiments that the above simplification and its associated solution scheme produce image reconstructions of the quality comparable to those obtained using dense sampling of the GPF phase.
\end{abstract}


\IEEEpeerreviewmaketitle

\section{Introduction}
\IEEEPARstart{T}{he} necessity to recover digital images from their distorted and noisy observations arises in a multitude of practical applications, with some specific examples including image denoising, super-resolution, image restoration, and watermarking, just to name a few \cite{30, 28, 29, 40}. In such cases, it is standard to assume that the observed image $v$ is formed by subjecting the original image $u$ to convolution with a point spread function\footnote{Note that, in optical imaging, this function is also referred to as an impulse transfer function \cite{39}.} (PSF) $i$, followed by contamination by white Gaussian noise (WGN) $\nu$. Thus, formally,
\begin{equation}\label{1}
\begin{split}
v=i \ast u + \nu.
\end{split}
\end{equation}
While $u$ and $v$ can be regarded as general members of the signal space $\mathbb{L}_2(\Omega)$ of real-valued functions on $\Omega \subseteq \mathbb{R}^2$, the PSF $i$ is normally a much smoother function, with effectively band-limited spectrum. As a result, the convolution with $i$ has a destructive effect on the informational content of $u$, in which case $v$ typically has a substantially reduced set of features with respect to $u$. This makes the problem of reconstruction of $u$ from $v$ a problem of significant practical importance \cite{41}.

Reconstruction of the original image $u$ from $v$ can be carried out within the framework of image deconvolution, which is a specific instance of a more general class of inverse problems \cite{42}. Most of such methods are Bayesian in nature, in which case the information lost in the process of convolution with $i$ is recovered by requiring the optimal solution to reside within a predefined functional class \cite{26, 27}. Thus, for example, in the case when $u$ is known to be an image of bounded variation, the above regularization leads to the famous Rudin-Osher-Fatemi reconstruction scheme, in which $u$ is estimated as a solution to the following problem \cite{0, 11}
\begin{equation} \label{2}
\hat{u} = \underset{u}{\arg\min} \left\{ \frac{1}{2} \| u \ast i - v \|_2^2 + \alpha \int |\nabla u| \, dx dy \right\},
\end{equation}
where $\alpha>0$ is the regularization parameter. It should be noted that, if the PSF obeys $\int i \, dx dy \neq 0$, the problem (\ref{2}) is strictly convex and therefore admits a unique minimizer, which can be computed by a spectrum of available algorithms \cite{0, 11}.

A particularly non-trivial version of deconvolution is commonly referred to as blind. In this case, the original image $u$ is to be estimated without the knowledge of the PSF \cite{42}. In this paper, however, we follow the philosophy of {\em hybrid deconvolution} \cite{Oleg07}, which takes advantage of any partial information on the PSF to improve the image reconstruction. Thus, in the algorithm described in this paper, the original image $u$ will be recovered from $v$ and some partial information on $i$.

Optical (and, in particular, turbulent) imaging is unarguably the field of applied sciences from which the notion of deconvolution has originally emanated \cite{Richardson72, Lucy74, 44}. In short-exposure imaging, however, computational methods of image restoration are still superseded by adaptive optics. As recently as a decade ago, the use of adaptive optics would have been considered as the only practical option. Nowadays, however, with the advent of distributed cluster computing and GPU-based image processing, it seems to be time to revisit the cost-to-performance characteristics of the existing tools of adaptive optics. Thus, in this work, our focus is on a specific tool of adaptive optics, known as the Shack-Hartmann interferometer \cite{12, 18}. Instead of completely excluding the interferometer from our measurement system, we propose to  modify its construction through reducing the number of its local wavefront lenses. Although the advantages of such a simplification are immediate to see, its main shortcoming is obvious as well: the smaller the number of lenses is, the stronger is the effect of undersampling and aliasing. Accordingly, to overcome this problem, we propose to augment the modified Shack-Hartmann interferometer by subjecting its output to the derivative compressed sensing (DCS) algorithm of \cite{14}. As it will be shown later in the paper, the PSF $i$ is determined by a generalized pupil function $P$, which can be expressed in a polar form as $P = A \, e^{\jmath \phi}$. While the amplitude $A$ can be measured via calibration or computed as a function of the aperture geometry, the phase $\phi$ is often influenced by environmental effects and hence it needs to be recovered from observations. It will be shown below that DCS is particularly well suited for reconstruction of $\phi$ from incomplete measurements of its partial differences. Such an estimate can be subsequently combined with $A$ to yield an estimate of the PSF $i$, which can in turn be used by a deconvolution algorithm. Thus, the proposed method for estimation of the PSF and subsequent deconvolution of $u$ can be regarded as a hybrid deconvolution technique, which comes to simplify the design and complexity of adaptive optics on the one hand, and to make the process of reconstruction of optical images as automatic as possible, on the other hand.

The rest of the paper is organized as follows. Section II summarizes basic technical preliminaries. In Section III, we describe the SH interferometer as well as phase measurements in optical imaging. In Section IV, we explain DCS and our new approach to solve it. In Section V we describe deconvolution process to recover the original image. Experimental results are presented in Section VI, while Section VII finalizes the paper with a discussion and main conclusions.

\section{Technical Preliminaries}
In short exposure imaging, due to aberrations in the imaging system induced by, e.g., atmospheric turbulence, the impulse response of an optical imaging system is often unknown \cite{15}. In order to better understand the setup under consideration, we first note that, in optical imaging, the PSF $i$ is obtained from an amplitude spread function (ASF) $h$ as $i := | h |^2$. The ASF, in turn, is defined in terms of the generalized pupil function (GPF) $P(x,y)$ as given by \cite{17}
\begin{equation} \label{3}
h(\xi,\eta)=\frac{1}{\lambda_w z_i}\int_{-\infty}^{\infty} \int_{-\infty}^{\infty} P(x,y) e^{-j\frac{2\pi}{\lambda z_i}( x \, \xi + y \, \eta)} \, dxdy,
\end{equation}
where $z_i$ is the focal distance and $\lambda_w$ is the optical wavelength. Being a complex-valued quantity, $P(x,y)$ can be represented in terms of its amplitude $A(x,y)$ and phase $\phi(x,y)$ as
\begin{equation} \label{4}
P(x,y) = A(x,y) \, e^{ \jmath \phi(x,y)}.
\end{equation}
Here, the GPF amplitude $A(x,y)$ (which is sometimes simply referred to as the aperture function) is normally a function of the aperture geometry. Thus, for instance, in the case of a circular aperture, $A(x,y)$ can be defined as \cite{15}
\begin{equation} \label{5}
A(r)=
 \begin{cases}
 1, \quad  &\mbox{if } r \le \frac{D}{2} \\
 0, \quad  &\text{otherwise}
 \end{cases}
\end{equation}
where $D$ denotes the pupil diameter. Thus, given $\phi(x,y)$, one could determine $h$ and therefore $i$. Unfortunately, the phase $\phi(x,y)$ does not have an analytic expression, and it has to be measured in practice using such tools as the Shack-Hartmann interferometer (SHI) \cite{12}.

As will be discussed later in the paper, the SHI is capable of sensing the partial derivatives of $\phi(x,y)$. Needless to say, in order to minimize the effect of aliasing on the estimation result, an accurate reconstruction of $\phi(x,y)$ requires taking a fairly large number of the samples of $\nabla \phi (x,y)$ \cite{Oleg2008}. In some applications, the number of sampling points (as defined by the number of local wavefront lenses) reaches as many as a few thousands. It goes without saying that reducing the number of lenses would have a positive impact on the SHI in terms of its cost and approachability. Alas, such a reduction is impossible without undersampling, which tends to have formidable effect on the overall quality of phase estimation.

In this paper, to minimize the effect of undersampling, we exploit DCS \cite{14}. As opposed to the classical compressed sensing (CCS) \cite{13}, in addition to the sparsifing constraints, DCS also uses constraints which are intrinsic in the definition of partial derivatives. Using these additional constraints -- which are called the cross-derivative constraints -- allows substantially improving the quality of reconstruction of $\phi(x,y)$, as compared to the case of CCS-based estimation.

\section{Shack-Hartmann Interferometer (SHI)}
As it was mentioned earlier, the SHI is typically used to measure the gradient $\nabla \phi(x,y)$ of the GPF phase $\phi(x,y)$, from which the values of the latter can be subsequently estimated. To this end, the unknown phase $\phi(x,y)$ is assumed to be expandable in terms of some basis functions $\{Z_k\}_{k=0}^\infty$, {\it viz.} \cite{18}
\begin{equation} \label{6}
\phi(x,y)=\sum_{k=0}^{\infty} a_k Z_k(x,y),
\end{equation}
where the representation coefficients $\{a_k\}_{k=0}^\infty$ are assumed to be unique and stably computable. Note that, in this case, the datum of $\{a_k\}_{k=0}^\infty$ uniquely identifies $\phi(x,y)$, while the coefficients $\{a_k\}_{k=0}^\infty$ can be estimated due to the linearity of (\ref{6}) which suggests
\begin{equation} \label{6.1}
\nabla \phi(x,y)=\sum_{k=0}^{\infty} a_k \, \nabla Z_k(x,y),
\end{equation}

The most frequent choice of $\{Z_k\}_{k=0}^\infty$ in adaptive optics is the Zernike polynomials (aka Zernike functions) \cite{17}. These polynomials constitute an orthonormal basis in the space of square-integrable functions defined over the unit disk in $\mathbb{R}^2$. Zernike polynomials can be subdivided in two subsets of the even $Z_n^m$ and odd $Z_n^{-m}$ Zernike polynomials which have very convenient analytical definitions as given by
\begin{align}
Z^{m}_n(\rho,\varphi) &= R^m_n(\rho)\,\cos(m\,\varphi) \\
Z^{-m}_n(\rho,\varphi) &= R^m_n(\rho)\,\sin(m\,\varphi)
\end{align}
where $m$ and $n$ are nonnegative integers with $n \ge m$, $0 \le \varphi < 2\pi$ is the azimuthal angle, and $0 \le \rho \le 1$ is the radial distance. The radial polynomials $R^m_n$ are defined as
\begin{equation}
R^m_n(\rho) = \! \sum_{k=0}^{(n-m)/2} \!\!\! \frac{(-1)^k\,(n-k)!}{k!\,((n+m)/2-k)!\,((n-m)/2-k)!} \;\rho^{n-2\,k}.
\end{equation}

Note that, since the Zernike polynomials above are defined using polar coordinates, it makes sense to re-express the phase $\phi$ and its gradient in the polar coordinate system as well. (Technically, this would amount to replacing $x$ and $y$ in (\ref{6})-(\ref{6.1}) by $\rho$ and $\varphi$, respectively.) Moreover, due to the property of the Zernike polynomials to be an orthonormal basis, the representation coefficients $\{a_k\}_{k=0}^\infty$ in in (\ref{6})-(\ref{6.1}) can be computed by orthogonal projection, namely
\begin{equation} \label{71}
a_k = \int_0^{2\pi} \int_0^1 \phi(\rho,\varphi) \, Z_k(\rho,\varphi) \, \rho \, d\rho \, d\varphi
\end{equation}
In practice, however, $\phi(\rho,\varphi)$ is unknown and therefore the coefficients $\{a_k\}_{k=0}^\infty$ need to be estimated by other means. Thus, in the case of the SHI, the coefficients can be estimated from a finite set of discrete measurements of $\nabla \phi(\rho,\varphi)$.

\begin{figure}[!t]
\begin{center}
\hspace{-5mm}\includegraphics[height=2.1in,width=2.1in]{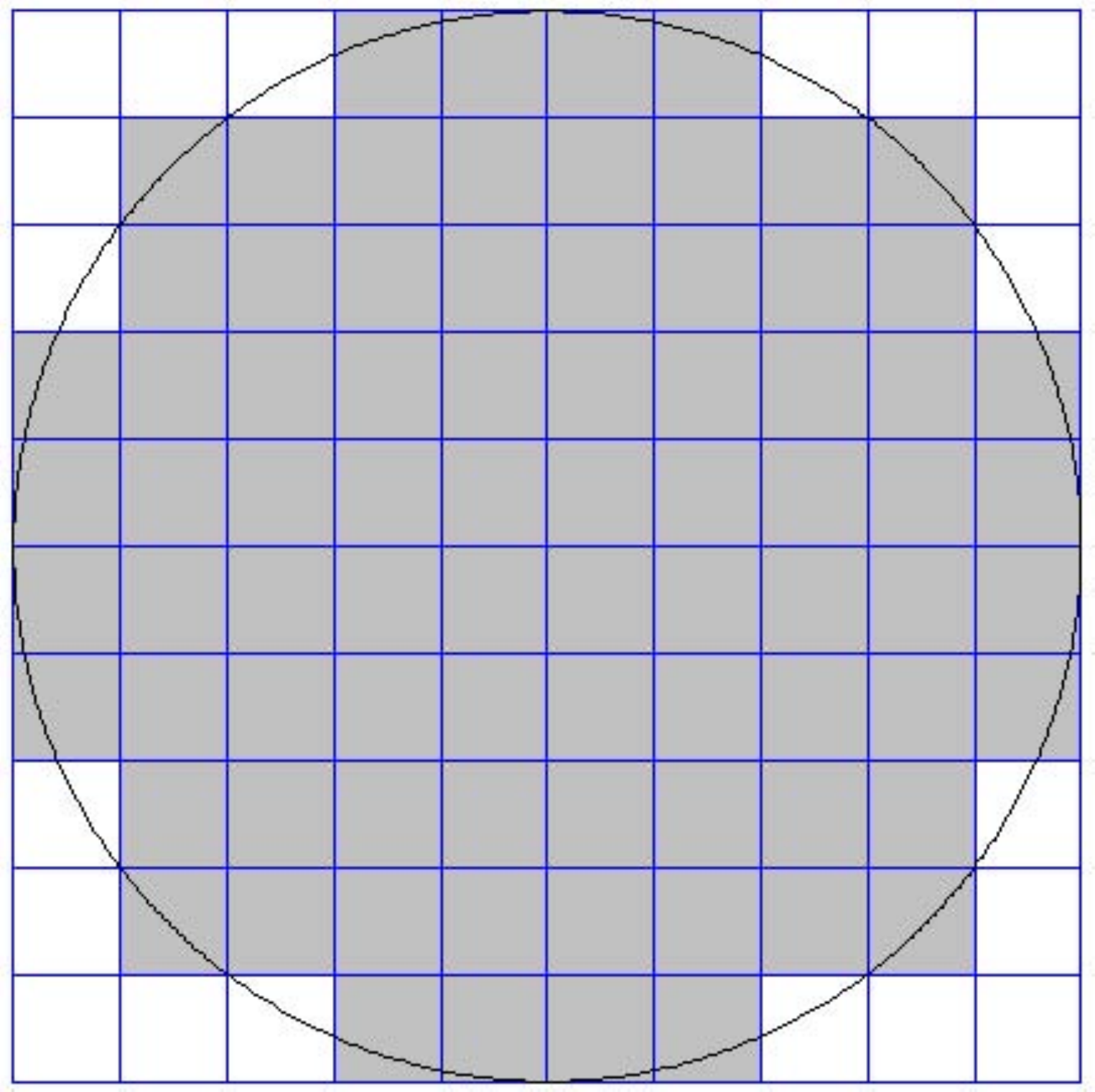}
\caption{An example of a $10 \times 10$ SHI array on a circular aperture. The shading indicates those blocks (i.e., lenses) which are rendered idle.} \label{F1}
\end{center}
\end{figure}

The main function of the SHI is to acquire discrete measurements of $\nabla \phi$ by means of linearization. The linearization takes advantage of subdividing a (circular) aperture into rectangular blocks with their sides formed by a uniform rectangular lattice. An example of such a subdivision is shown in Fig.~\ref{F1} for the case of a $10 \times 10$ lattice grid. Subsequently, it is assumed that the grid is sufficiently fine to approximate a restriction of the phase $\phi$ to each of the above blocks by a linear function. This results in a piecewise linear approximation of $\phi$, whose accuracy improves asymptotically when the lattice size goes to infinity\footnote{More rigorously, one can show that, as long as $\phi$ is uniformly continuous, its piecewise linear approximation converges uniformly, as the grid size goes to infinity.}. Formally, let $\Omega := \{(x,y) \in \mathbb{R}^2 \mid x^2+y^2 \le D^2\}$ be a circular aperture of radius $D$ and $\mathcal{S} = \{(x,y) \in \mathbb{R}^2 \mid \max\{|x|,|y|\} \le D \}$ be a square subset of $\mathbb{R}^2$ such that $\Omega \subset \mathcal{S}$. Then, for each polar coordinate $(\rho,\varphi) \in \Omega$ and an $N \times N$ grid of square blocks of size $2 D / N \times 2 D / N$, the phase $\phi$ can be expressed as
\begin{equation} \label{8}
\phi(x,y) \approx ax + by + c,
\end{equation}
for all $(x,y)$ in a neighbourhood of $(\rho \, \cos \varphi, \rho \, \sin \varphi)$. The approximation in (\ref{8}) suggests that
\begin{equation} \label{9}
\begin{split}
\nabla\phi(x,y) \approx [a, b]^T
\end{split}
\end{equation}
where $(\cdot)^T$ denotes matrix transposition. While $c$ in (\ref{8}) can be derived from boundary conditions, coefficients $a$ and $b$ should be determined via direct measurements. To this end, the SHI is endowed with an array of small focusing lenses, which are supported over each of the square blocks of the discrete grid. In the absence of phase aberrations, the focal points of the lenses are spatially identified and registered using a high-resolution CCD detector, whose imaging plane is aligned with the plane of the focal points. Then, when the wavefront gets distorted as a result of, e.g., atmospheric turbulence, the focal points are ``pushed" towards new spatial positions, which can also be pinpointed by the same detector. The resulting displacements can therefore be measured and subsequently related to the values of $\nabla \phi$ at corresponding points.

\begin{figure}[!t]
\centering
\includegraphics[width = 7cm, height = 5cm]{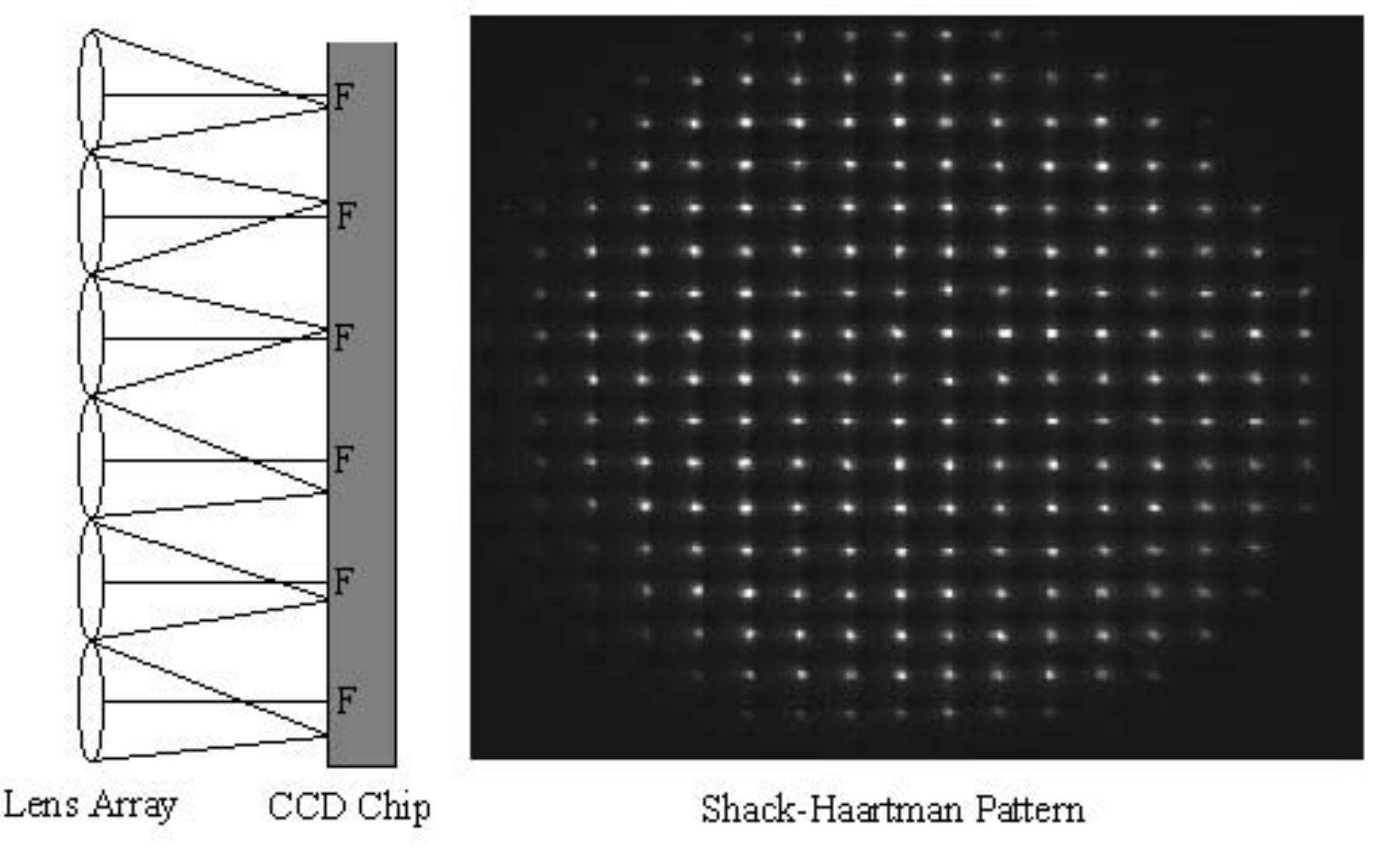}
\label{fig:tabsubfig21}
\caption{Basic structure of the SHI and a resulting pattern of the focal points.}
\label{fig222}
\end{figure}

To explain how the above procedure can be performed, additional notations are in order. Let $\Omega_d$ denote a finite set of spatial coordinates defined according to
\begin{align}
\Omega_d &:= \Big\{ (x_d,y_d) \in \Omega \,\, \big| \notag\\
x_d &=-D+\frac{2D}{N} \left(i+\frac{1}{2}\right), \,\, i=0,1,\ldots,N-1 \\
y_d &=-D+\frac{2D}{N} \left(j+\frac{1}{2}\right), \,\, j=0,1,\ldots,N-1 \notag \\
&\mbox{and } x_d^2 + y_d^2 \le D^2 \Big\}. \notag
\end{align}
Note that the set $\Omega_d$ can be thought of as a set of the spatial coordinates of the geometric centres of the wavefront lenses, restricted to the domain of aperture $\Omega$. Under some reasonable assumptions, one can then show \cite{19} that the focal displacement $\Delta(x,y) = [\Delta_x(x,y), \Delta_y(x,y)]^T$ measured at some $(x,y) \in \Omega_d$ is related to the value of $\nabla \phi(x,y)$ as given by
\begin{equation}
\nabla \phi(x,y) \approx \frac{1}{f} \Delta \phi(x,y), \quad \forall (x,y) \in \Omega_d,
\end{equation}
where $f$ is the focal distance of the wavefront lenses. Such a measurement setup is depicted in Fig.~\ref{fig222} along with an example of measured focal points.

Now, given a total of $M : = \# \Omega_d$ measurements of $\nabla \phi$ over $\Omega_d$, one can try to recover a useful approximation of $\phi$ over the whole $\Omega$ in the form of a projection of $\phi$ on the linear subspace spanned by all Zernike polynomials up to the order $L$ inclusive. In this case, it is possible to estimate the representation coefficients $\{a_k\}_{k=0}^L$ via solution of
\begin{align} \label{E1}
\min_{\{a_k\}} \sum_{(x,y) \in \Omega_d} \big \| \sum_{k=0}^L a_k \nabla Z_k(x,y) - f^{-1} \Delta (x,y) \big\|_2^2
\end{align}
subject to appropriate boundary conditions. It is worthwhile noting that (\ref{E1}) can be rewritten in a vector-matrix form as
\begin{equation} \label{E2}
\min_{\bf a} \| {\bf Z} \, {\bf a}  - {\bf d} \|_2^2, \quad \mbox{s.t.} \,\, {\bf a} \succeq 0,
\end{equation}
where $\bf Z$ is a $2 M \times L+1$ matrix composed of the values of the partial derivatives of the Zernike polynomials, $\bf d$ is a measurement (column) vector of length $2 M$, and ${\bf a} = [a_0, a_1, \dots, a_L]^T$ is a vector of the representation coefficients of $\phi$. The constraint ${\bf a} \succeq 0$ in (\ref{E2}) is supposed to further regularize the solution by forcing $\bf a$ to be  a member of some convex cone as well. Thus, for example, if the mean value of $\phi$ can be assumed to be equal to zero, the solution to (\ref{E2}) can be computed as
\begin{equation}
{\bf a} = {\bf Z}^{\#} {\bf d},
\end{equation}
where ${\bf Z}^{\#}$ denotes the pseudo-inverse of $\bf Z$, whose definition is unique and stable as long as the row-rank of $\bf Z$ is greater or equal to $L+1$ (hence suggesting that $2 M \ge L+1$). Having estimated $\bf a$, the phase $\phi$ can be approximated as
\begin{equation}
\phi(\rho, \varphi) \approx \sum_{k=0}^L a_k Z_k(\rho, \varphi).
\end{equation}

A higher accuracy of phase estimation requires using higher-order Zernike polynomials, which in turn necessitates a proportional increase in the number of wavefront lenses. Moreover, as required by the linearization procedure in the SHI, the lenses have to be of a relatively small sizes (sometimes, on the order of a few microns), which may lead to the use of a few thousand lenses per one interferometer. Accordingly, to simplify the construction and to reduce the cost of SHIs, we propose to substantially reduce the number of wavefront lenses, while compensating for the induced information loss through the use of DCS, which is detailed next.

\section{Derivative Compressive Sampling}
\subsection{Classical Compressed Sensing}
Central to signal processing is the Shannon-Nyquist theorem \cite{45}, which specifies conditions on which a band-limited signal can be stably and uniquely recovered from its discrete measurements. However, in around 2005, a different sampling theorem was formed which, in some cases, abrogates the fundamentals of its predecessor. This new theory, nowadays known as compressed sensing (aka compressive sampling), asserts that signals, which admit a sparse representation in a predefined basis/frame, can be recovered from their discrete measurements, whose number is proportional to the $\ell_0$-norm of the coefficients of the sparse representation. In such a case, the sparser the representation of the signal is, the smaller can be the number of  measurements required for signal reconstruction. As a result, cases are numerous in which the conditions of compressed sensing far superseded those of the Shannon-Nyquist sampling \cite{4, 13}.

Despite its widespread success in countless applications, the theory of compressed sensing is not entirely free of limitations. One of such limitations stems from the necessity to use a non-linear decoder. Indeed, while in the case of classical sampling, the reconstruction of a time-domain signal is implemented through linear interpolation (i.e., linear filtering), in the case of compressed sensing, the reconstruction involves solution of a convex optimization problem. Specifically, let us consider a typical setup of {\em classical} compressed sensing (CCS), in which $\yy \in \mathbb{R}^m$ represents an observed version of $\xx \in \mathbb{R}^n$, related according to
\begin{equation} \label{13}
\yy = \Psi \xx,
\end{equation}
where $\Psi \in \mathbb{R}^{m\times n}$ is an observation (sampling) matrix with $n > m$.

The recovery of $\xx$ from $\yy$ based on (\ref{13}) is impossible to implement in a unique and stable way, unless it is known that $\xx$ is sparse and hence has a relatively low value of $\|\xx\|_0$. In such a case, if the sampling matrix $\Psi$ satisfies the restricted isometry property (RIP) \cite{4, 13} with respect to a certain class of sparse signals to which $\xx$ is assumed to belong, then CCS recovers $\xx$ as a solution to \cite{20, 21}
\begin{equation}
\xx = \arg \min_{\xx^\prime} \left\{ \| \xx^\prime \|_1 \mid \Psi \xx^\prime = \yy \right\},
\end{equation}
which is a convex minimization problem, which is straightforward to reformulate in terms of linear programming. Moreover, in the case when the measurements $\yy$ are error-prone, a more robust version of CCS is to recover $\xx$ as given by
\begin{equation}\label{E3}
\xx = \arg \min_{\xx^\prime} \left\{ \| \xx^\prime \|_1 \mid \| \Psi \xx^\prime - \yy \|_2^2 \le \epsilon \right\},
\end{equation}
where $\epsilon > 0$ is a parameter controlling the size of the noise. Moreover, it was shown in \cite{4, 13}, that the estimation error in the signal reconstructed according to (\ref{E3}) can be bounded by a linear function of $\epsilon$. This implies robustness of the CCS reconstruction towards the presence of measurement noise.

It should be finally noted that the optimization problem (\ref{E3}) can be reformulated in its alternative Lagrangian form, in which case one can find
\begin{equation} \label{17}
\xx = \arg \min_{\xx^\prime} \left\{ \frac{1}{2} \| \Psi \xx^\prime - \yy \|_2^2 + \lambda \| \xx^\prime \|_1 \right\}
\end{equation}
where $\lambda > 0 $ is an optimal Lagrange multiplier \cite{20}. In what follows, it is assumed that an optimal value of $\lambda$ is known. (For more details on this subject, the reader is referred to \cite{20} as well as to the later sections of this paper).

\subsection{Derivative Compressed Sensing (DCS)}
Let the partial derivatives of $\phi$ evaluated (by means of the SHI) at the points of set $\Omega_d$ be column-stacked into vectors $\ffx$ and $\ffy$ of length $M = \#\Omega_d$.  In what follows, the partial derivatives $\ffx$ and $\ffy$ are assumed to be sparsely representable by an orthonormal basis in $\mathbb{R}^M$. Representing such a basis by an $M \times M$ unitary matrix $W$, the above assumption suggests the existence of two {\em sparse} vectors $\ccx$ and $\ccy$ such that $W \ccx$ and $W \ccy$ provide accurate approximations of the partial derivatives of the original phase $\phi$ evaluated over $\Omega_d$.

Now, the simplification of the SHI proposed in the current paper amounts to reducing the number of wavefront lenses to a minimum. Formally, such a reduction can be described by two $n \times M$ sub-sampling matrices $\Psi_x$ and $\Psi_y$, where $n < M$. Specifically, let $\bbx := \Psi_x \ffx$ and $\bby := \Psi_y \ffy$ be incomplete (partial) observations of $\ffx$ and $\ffy$, respectively. Then, the noise-free counterparts of the partial derivatives can be approximated by $W \ccx^\ast$ and $W \ccy^\ast$, respectively, where
\begin{equation}
\ccx^\ast = \arg \min_{\ccx^\prime} \left\{ \frac{1}{2} \| \Psi_x W \ccx^\prime - \bbx \|_2^2  + \lambda_x \| \ccx^\prime \|_1 \right\}
\end{equation}
and
\begin{equation}
\ccy^\ast = \arg \min_{\ccy^\prime} \left\{ \frac{1}{2} \| \Psi_y W \ccy^\prime - \bby \|_2^2  + \lambda_y \| \ccy^\prime \|_1 \right\}
\end{equation}
for some $\lambda_x, \lambda_y > 0$. Moreover, in the case when $\lambda_x = \lambda_y$, the above estimates can be combined together. To this end, let $\cc = [\ccx, \ccy]^T$, $\bb = [\bbx, \bby]^T$, and $A = \mbox{diag} \{\Psi_x W, \Psi_y W\} \in \mathbb{R}^{2 n \times 2 M}$. Then,
\begin{equation} \label{E4}
\cc^\ast = \arg \min_{\cc^\prime} \left\{ \frac{1}{2} \| A \cc^\prime - \bb \|_2^2  + \lambda \| \cc^\prime \|_1 \right\},
\end{equation}
where $\lambda=\lambda_x = \lambda_y$. In this form, the problem (\ref{E4}) is identical to (\ref{17}) and hence it can be solved by a variety of optimization algorithms \cite{20, 21}.

The DCS algorithm augments CCS by subjecting the minimization in (\ref{E4}) to an additional constraint which stems from the fact that \cite{14}
\begin{equation} \label{20}
\frac{\partial^2 \phi}{\partial x \, \partial y}=\frac{\partial^2 \phi}{\partial y \, \partial x},
\end{equation}
which is valid for any two times continuously differentiable $\phi(x,y)$. Thus, in particular, the constraint implies the existence of two partial differences matrices $D_x$ and $D_y$ which obey
\begin{equation} \label{E5}
D_x \ffy = D_y \ffx.
\end{equation}
Consequently, if $T_x$ and $T_y$ are the matrices satisfying $T_x \cc = \ccx$ and $T_y \cc = \ccy$, respectively, then (\ref{E5}) can be re-expressed in terms of $\cc$ as
\begin{equation}
D_y T_x \cc = D_x T_y \cc
\end{equation}
or
\begin{equation}
B \cc = 0,
\end{equation}
where $B := D_y T_x - D_x T_y$. Thus, DCS solves the constrained minimization problem as given by
\begin{align} \label{E5}
\cc^\ast = \arg \min_{\cc^\prime} &\left\{ \frac{1}{2} \| A \cc^\prime - \bb \|_2^2  + \lambda \| \cc^\prime \|_1 \right\}, \\
&\mbox{s.t.  } B \cc^\prime = 0 \notag
\end{align}

A solution to (\ref{E5}) can be found, for instance, by means of the Bregman algorithm \cite{22}, in which case $\cc^\ast$ is approximated by a stationary point of the sequence of iterations produced by
\begin{equation} \label{E6}
\begin{cases}
\cc^{(t+1)} = \arg \min_{\cc^\prime} \Big\{ \frac{1}{2} \| A \cc^\prime - \bb \|_2^2  + \\
\hspace{2.5cm} + \lambda \| \cc^\prime \|_1 + \frac{\delta}{2} \| B \cc^\prime + p^{(t)} \|_2^2 \Big\} \\
p^{(t+1)} = p^{(t)} + \delta B \cc^{(t+1)},
\end{cases}
\end{equation}
where $p^{(t)}$ is a vector of Bregman variables (or, equivalently, augmented Lagrange multipliers) and $\delta > 0$ is a user-defined parameter\footnote{In this work, we use $\delta = 0.5$.}.

The $\cc$-update step in (\ref{E6}) has the format of a standard basis pursuit de-noising (BPDN) problem \cite{46}, which can be solved by a variety of optimization methods \cite{23}. In the present paper, we used the FISTA algorithm of \cite{33} due to the simplicity of its implementation as well as for its remarkable convergence properties. It should be noted that the algorithm does not require explicitly defining the matrices $A$ and $B$. Only the {\em operations} of multiplication by these matrices and their transposes need to be known, which can be implemented in an implicit and computationally efficient manner.

Once an optimal $\cc^\ast$ is recovered, it can be used to estimate the noise-free versions of $\ffx$ and $\ffy$ as $W T_x \cc^\ast$ and $W T_y \cc^\ast$, respectively. These estimates can be subsequently passed on to the fitting procedure of Section III to recover the values of $\phi$, which, in combination with a known aperture function $A$, provide an estimate of the PSF $i$ as an inverse discrete Fourier transform of the autocorrelation of $P = A \, e^{\jmath \phi}$. Algorithm 1 below summarizes our method of estimation of the PSF.

\begin{algorithm}
\setlength{\leftmargini}{0pt}
\caption{PSF estimation via DCS}
\begin{enumerate}
\item {\it Data:} $\bbx$, $\bby$, and $\lambda > 0$
\item {\it Initialization:} For a given transform matrix $W$ and matrices/operators $\Psi_x$, $\Psi_y$, $D_x$, $D_y$, $T_x$ and $T_y$, preset the procedures of multiplication by $A$, $A^T$, $B$ and $B^T$.
\item {\it Phase recovery:} Starting with an arbitrary $\cc^{(0)}$ and $p^{(0)} = 0$, iterate (\ref{E6}) until convergence to result in an optimal $\cc^\ast$. Use the estimated (full) partial derivatives $W T_x \cc^\ast$ and $W T_y \cc^\ast$ to recover the values of $\phi$ over $\Omega$.
\item {\it PSF estimation:} Using a known aperture function $A$, compute the inverse Fourier transform of $P = A \, e^{\jmath \phi}$ to result in a corresponding ASF $h$. Estimate the PSF $i$ as $i = |h|^2$.
\end{enumerate}
\label{algo1}
\end{algorithm}

The estimated PSF can be used to recover the original image $u$ from $v$ through the process of deconvolution as explained in the section that follows.

\section{Deconvolution}
The acquisition model \eqref{1} can be rewritten in an equivalent operator form as given by
\begin{equation} \label{30}
v = \mathcal{H}\{u\}+\nu,
\end{equation}
where $\mathcal{H}$ denote the operator of convolution with the estimated PSF $i$. Note that, in this case, the noise term $\nu$ accounts for both measurement noise as well as the inaccuracies related to estimation error in $i$.

The deconvolution problem of finding a useful approximation of $u$ given its distorted measurement $v$ can be addressed in many way, using a multitude of different techniques \cite{31, 32, 33}. In this work, we use the ROF model and recover a regularized approximation as a solution of
\begin{equation} \label{E7}
u^\ast = \arg \min_u \left\{ \frac{1}{2} \| \mathcal{H}\{u\} -v \|_2^2 + \gamma \, \| u \|_{TV} \right\},
\end{equation}
where $\| u \|_{TV} = \int \int | \nabla  u | \, dx \, dy $ denotes the total variation (TV) semi-norm of $u$.

One computationally efficient way to solve \eqref{E7} is to substitute a direct minimization of the cost function in \eqref{E7} by recursively minimizing a sequence of its local quadratic majorizers \cite{33}. In this case, the optimal solution $u^\ast$ can be approximated by the stationary point of a sequence of intermediate solutions produced by
\begin{equation} \label{E8}
\begin{cases}
w^{(t)} = u^{(t)} + \mu \, \mathcal{H}^\ast \left\{ v - \mathcal{H}\{u^{(t)}\}  \right\} \\
u^{(t+1)} = \arg \min_u \left\{ \frac{1}{2} \| u - w^{(t)} \|_2^2 + \gamma \, \| u \|_{TV} \right\},
\end{cases}
\end{equation}
where $\mathcal{H}^\ast$ is the adjoint of $\mathcal{H}$ and $\mu$ is chosen to satisfy $\mu > \| \mathcal{H}^\ast \mathcal{H} \|$. In this paper, the TV denoising at the second step of \eqref{E8} has been performed using the fixed-point algorithm of Chambolle \cite{11}. The convergence of \eqref{E8} can be further improved by using the same FISTA algorithm of \cite{33}. The resulting procedure is summarized below in Algorithm 2.

\begin{algorithm}
\setlength{\leftmargini}{0pt}
\caption{TV deconvolution using FISTA}
\begin{enumerate}
\item {\it Initialize:} Select an initial value $u^{(0)}$; set $y^{(0)} = u^{(0)}$ and $\tau^{(0)}=1$\\
\item {\it Repeat until convergence:}
\begin{itemize}
\item $w^{(t)} = y^{(t)} + \mu \, \mathcal{H}^\ast \left\{ v - \mathcal{H}\{y^{(t)}\}  \right\}$
\item $u^{(t+1)} = \arg \min_u \left\{ \frac{1}{2} \| u - w^{(t)} \|_2^2 + \gamma \, \| u \|_{TV} \right\}$
\item $\tau^{(t+1)} = 0.5  \left( 1+\sqrt{1+ 4 \, (\tau^{(t)})^2} \right)$
\item $y^{(t+1)} = u^{(t+1)} + (\tau^{(t)}/\tau^{(t+1)}) (u^{(t+1)} - u^{(t)})$
\end{itemize}
\end{enumerate}
\label{algo2}
\end{algorithm}

In summary, Algorithms 1 and 2 represent the essence of the proposed algorithm for hybrid deconvolution of optical images. The next section provides experimental results which further support the value and applicability of the proposed methodology.

\section{Results}
To demonstrate the viability of the proposed approach, its performance has been compared against reference methods. The first reference method used a dense sampling of the phase (as it would have been the case with a regular SHI), thereby eliminating the need for a CS-based phase reconstruction. The resulting method is referred below to as the dense sampling (DS) approach. Second, to assess the importance of incorporation of the cross-derivative constraints, we have used both CCS and DCS for phase recovery. In what follows, comparative results for phase estimation and subsequent deconvolution are provided for all the above methods.

\subsection{Phase recovery}
\begin{figure}[t]
\centering
\subfigure[]
{
\includegraphics[width = 4cm]{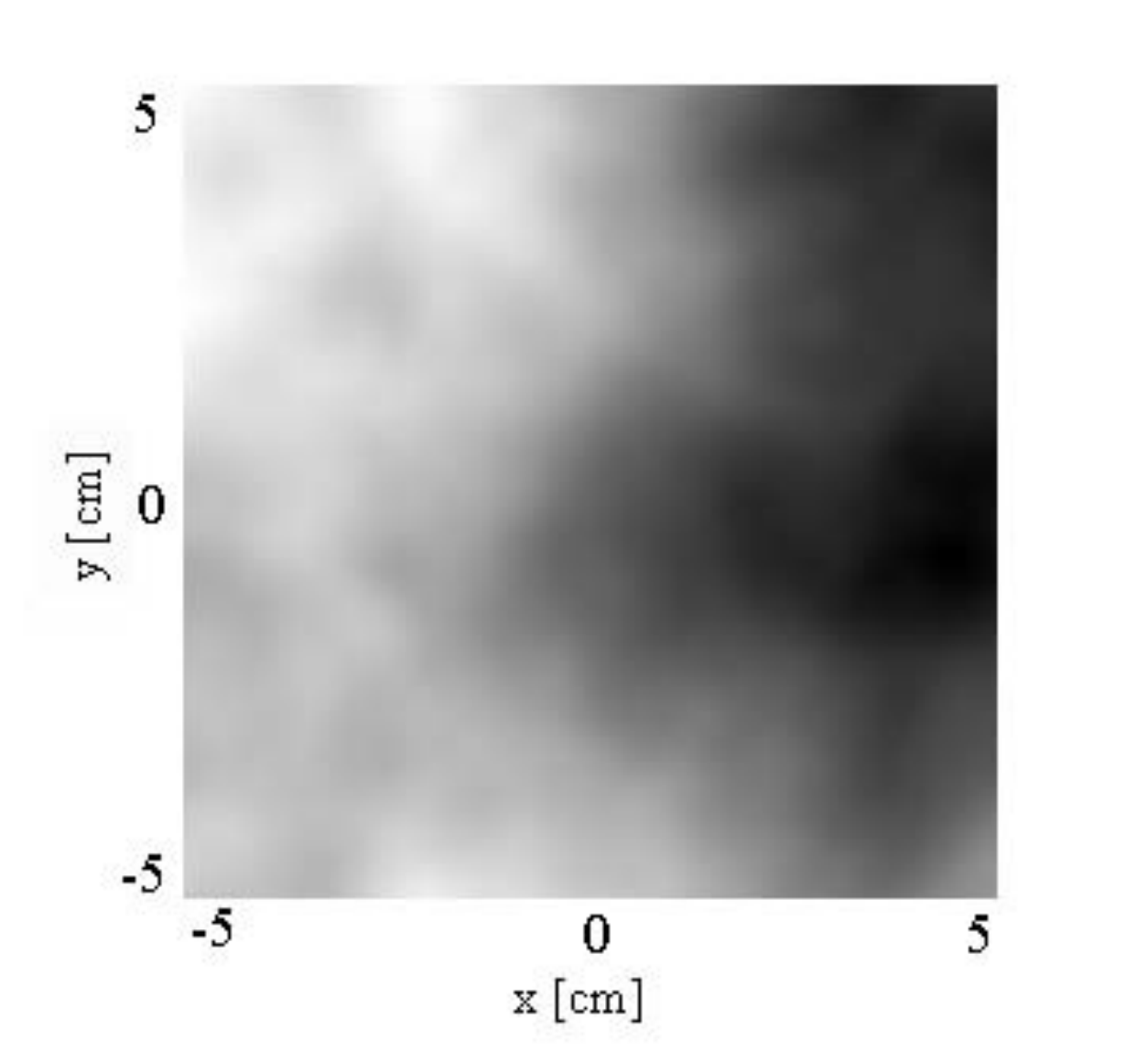}
\label{fig:tabsubfig1}
}
\\
\subfigure[]{
\includegraphics[width = 4cm]{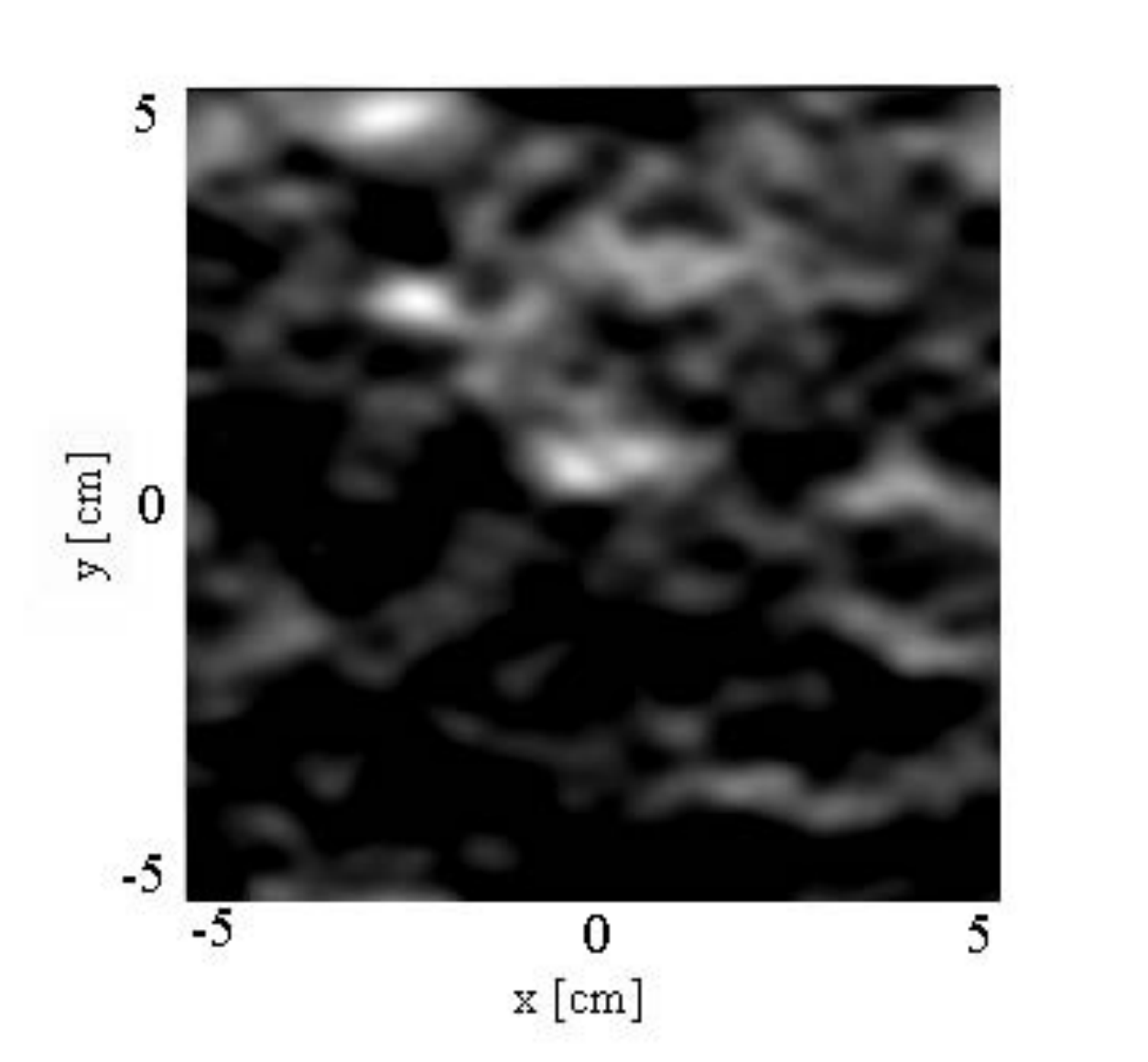}
\label{fig:tabsubfig2}
}
\subfigure[]{
\includegraphics[width = 4cm]{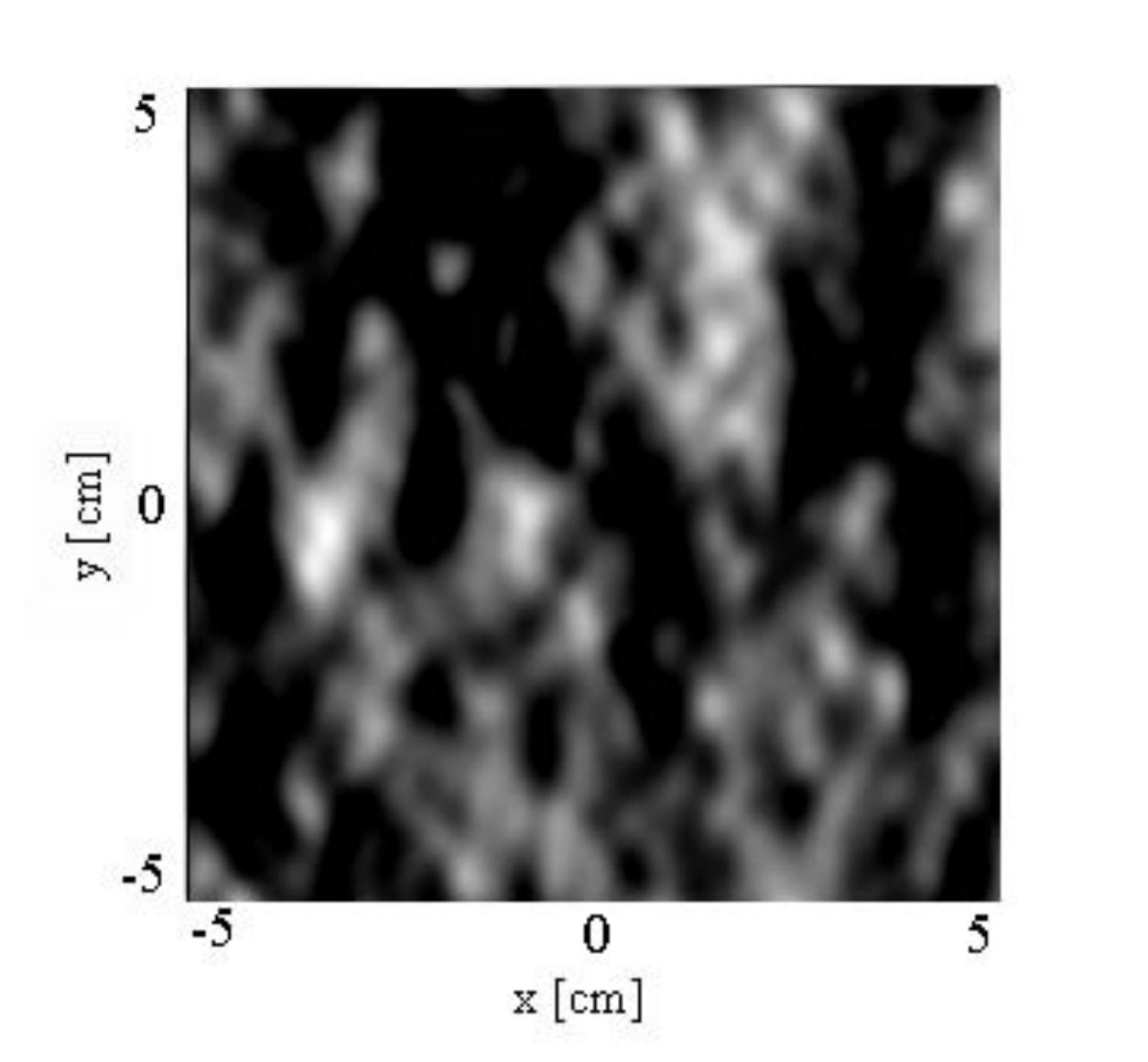}
\label{fig1:tabsubfig3}
}
\caption{An example of a simulated phase $\phi$ (subplot (a)) along with its partial derivatives w.r.t. $x$ (subplot (b)) and $y$ (subplot (c)).}
\label{F3}
\end{figure}

To assess the performance of the proposed and reference methods under controllable conditions, simulation data was used. The random nature of atmospheric turbulence necessitated the use of statistical methods to model its effect on a wavefront propagation. Specifically, in this paper, the effect of atmospheric turbulence has been described using the modified Von Karman PSF model\cite{35}. A typical example of a GPF phase $\phi$ is shown in subplot (a) of Fig.~\ref{F3}. In the shown case, the size of the phase screen was set to be equal to $10 \times 10$ cm, while the sampling was performed over a $128 \times 128$ uniform lattice (which would have corresponded to the use of 16384 lenses of a SHI). The partial derivatives $\partial \phi / \partial x$ and $\partial \phi / \partial y$ are shown in subplots (b) and (c) of Fig.~\ref{F3}, respectively.

In the present paper, the subsampling matrices $\Psi_x$ and $\Psi_y$ were obtained from an identity matrix $I$ through a random subsampling of its rows to result in a required compression ratio $r$ (to be specified below). To sparsely represent the partial derivatives of $\phi$, $W$ was defined to correspond to a four-level orthogonal wavelet transform using the nearly symmetric wavelets of Daubechies with five vanishing moments \cite{92}.

\begin{figure}[!t]
\centering
\includegraphics[width = 9cm, height = 6cm]{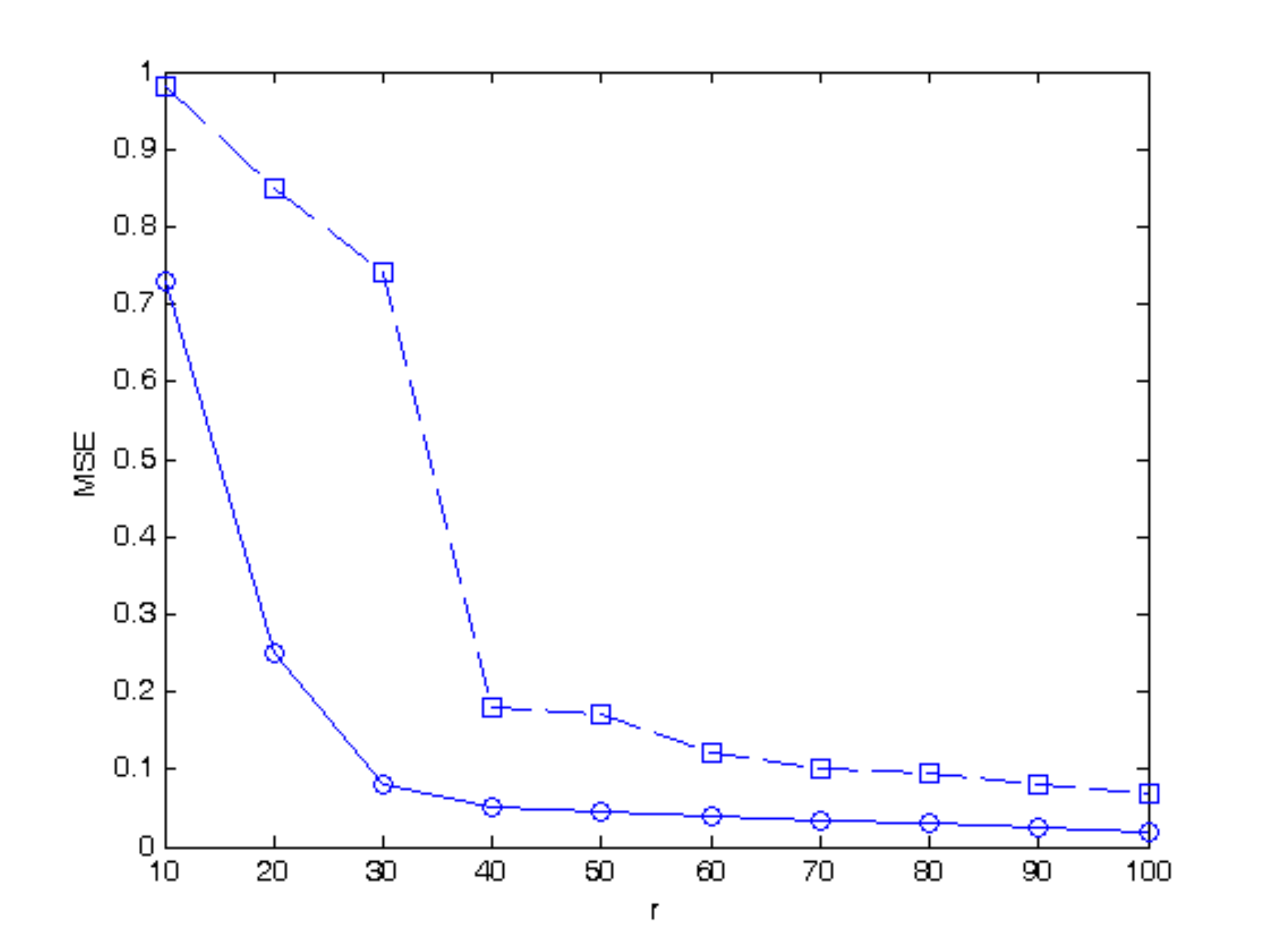}
\label{fig1:tabsubfig2}
\caption{MSE of phase reconstruction obtained with different methods as a function of $r$. Here, the dashed and solid lines correspond to CCS and DCS, respectively, and SNR is equal to 40 dB.}
\label{F4}
\end{figure}

To demonstrate the value of using the cross derivative constraint for phase reconstruction, the CCS and DCS algorithms have been compared in terms of the mean square errors (MSE) of their corresponding phase estimates. The results of this comparison are summarized in Fig.~\ref{F4} for different compression ratios (or, equivalently, (sub)sampling densities) and SNR = 40 dB.

As expected, one can see that DCS results in lower values of MSE as compared to CCS, which implies a higher accuracy of phase reconstruction. Moreover, the difference in the performances of  CCS and DCS appears to be more significant for lower sampling rates, while both algorithms tend to perform similarly when the sampling density approaches the DS case. Specifically, when the sampling density is equal to $r=0.3$, DCS results in a ten times smaller MSE as compared to the case of CCS, whereas both algorithms have comparable performance for $r = 0.83$. This result characterizes DCS as a better performer than CCS in the case of relatively low (sub)sampling rates.

\begin{figure}[!t]
\subfigure[]{
\includegraphics[width = 4cm]{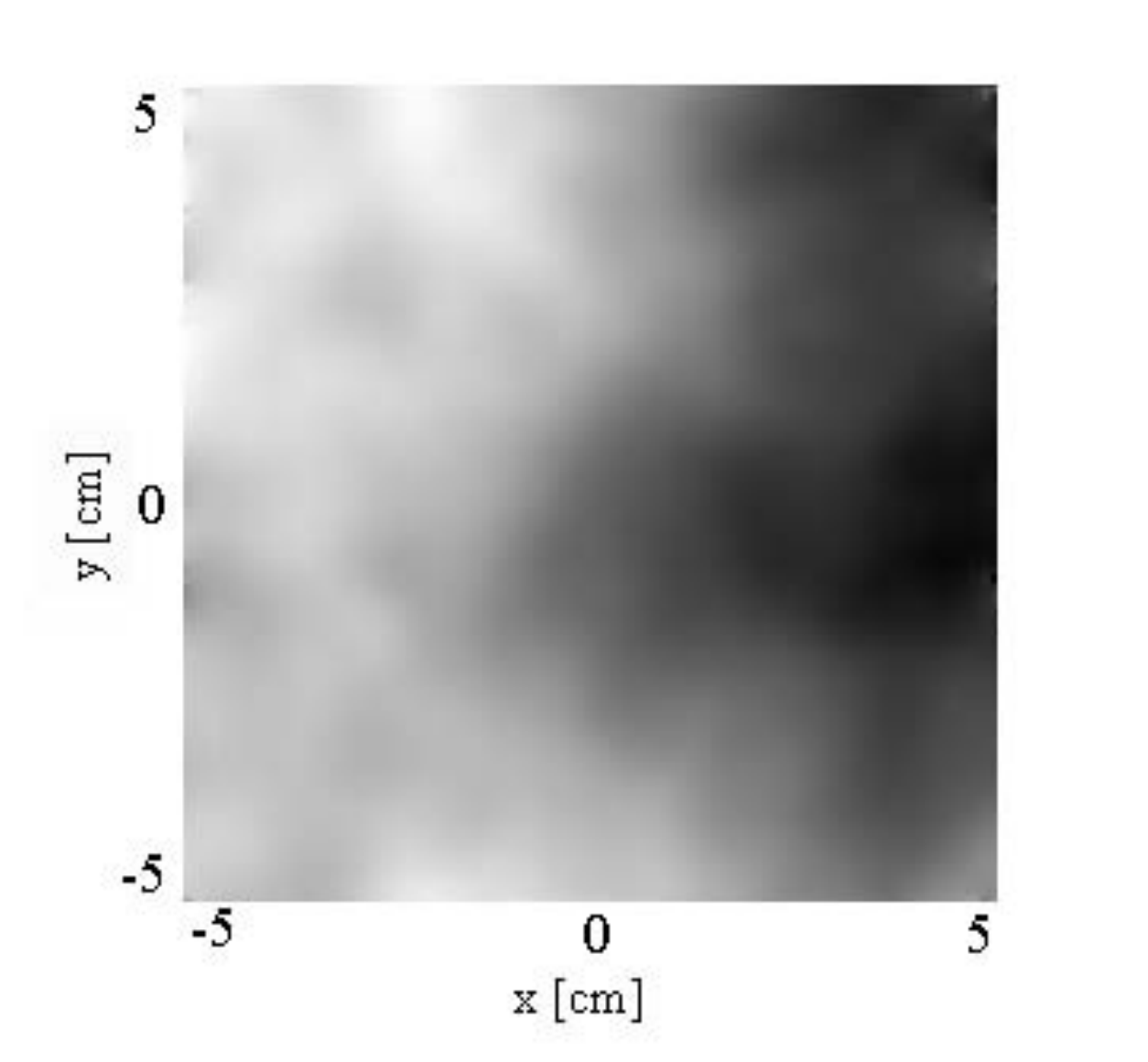}
\label{fig:tabsubfig21}
}
\subfigure[]{
\includegraphics[width = 4cm]{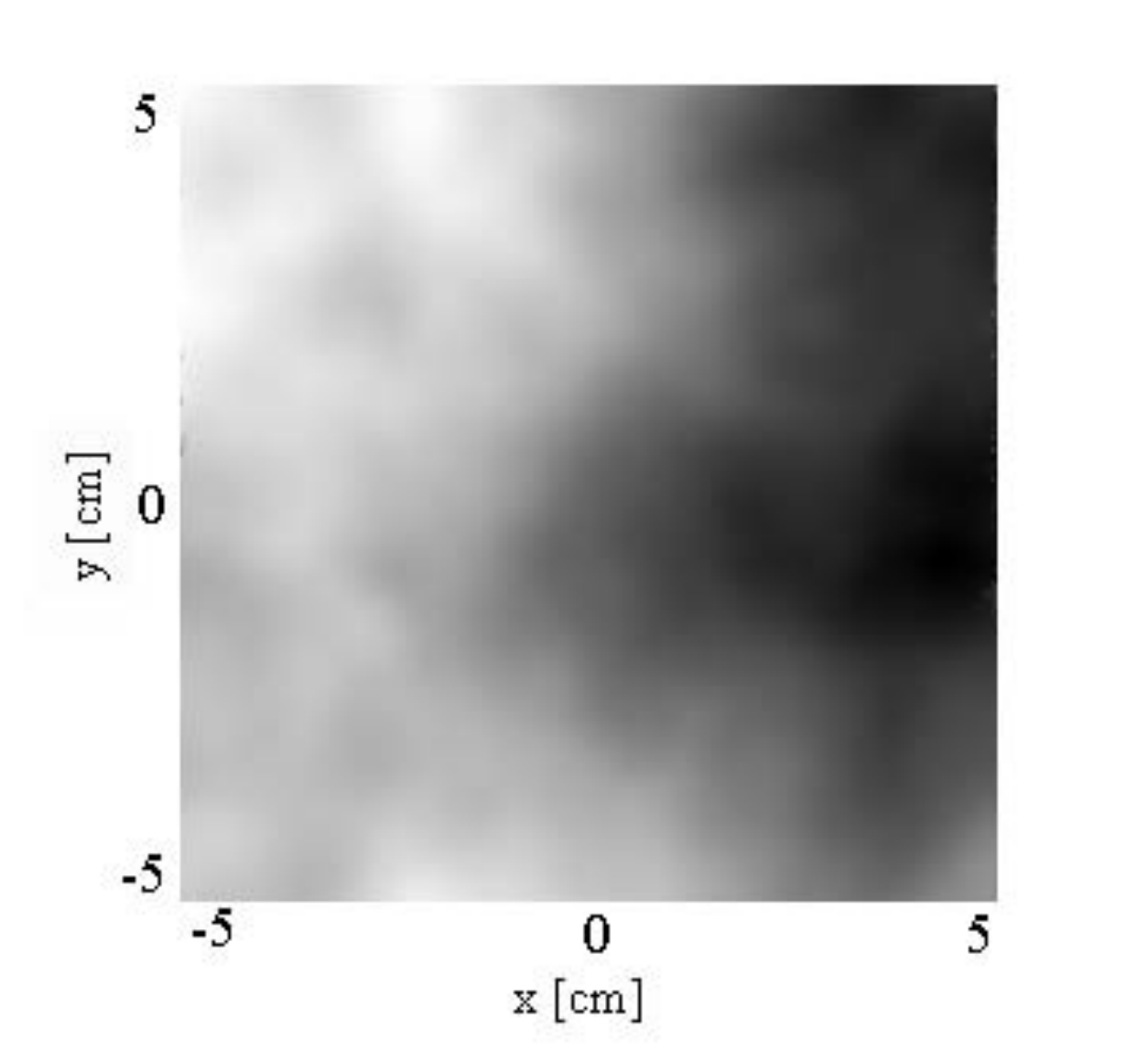}
\label{fig:tabsubfig31}
}\\
\subfigure[]{
\includegraphics[width = 4cm]{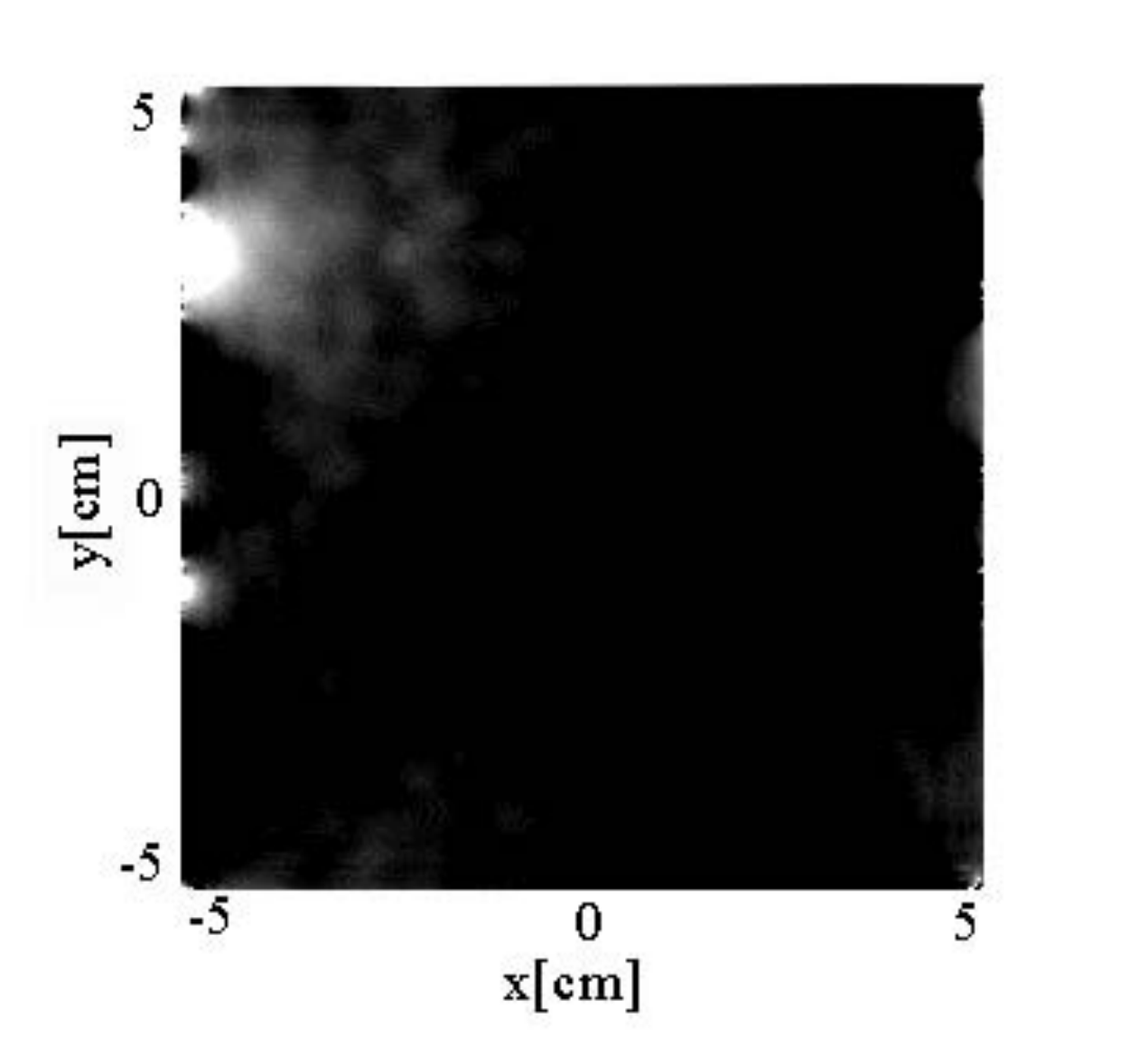}
\label{fig:tabsubfig31}
}
\subfigure[]{
\includegraphics[width = 4cm]{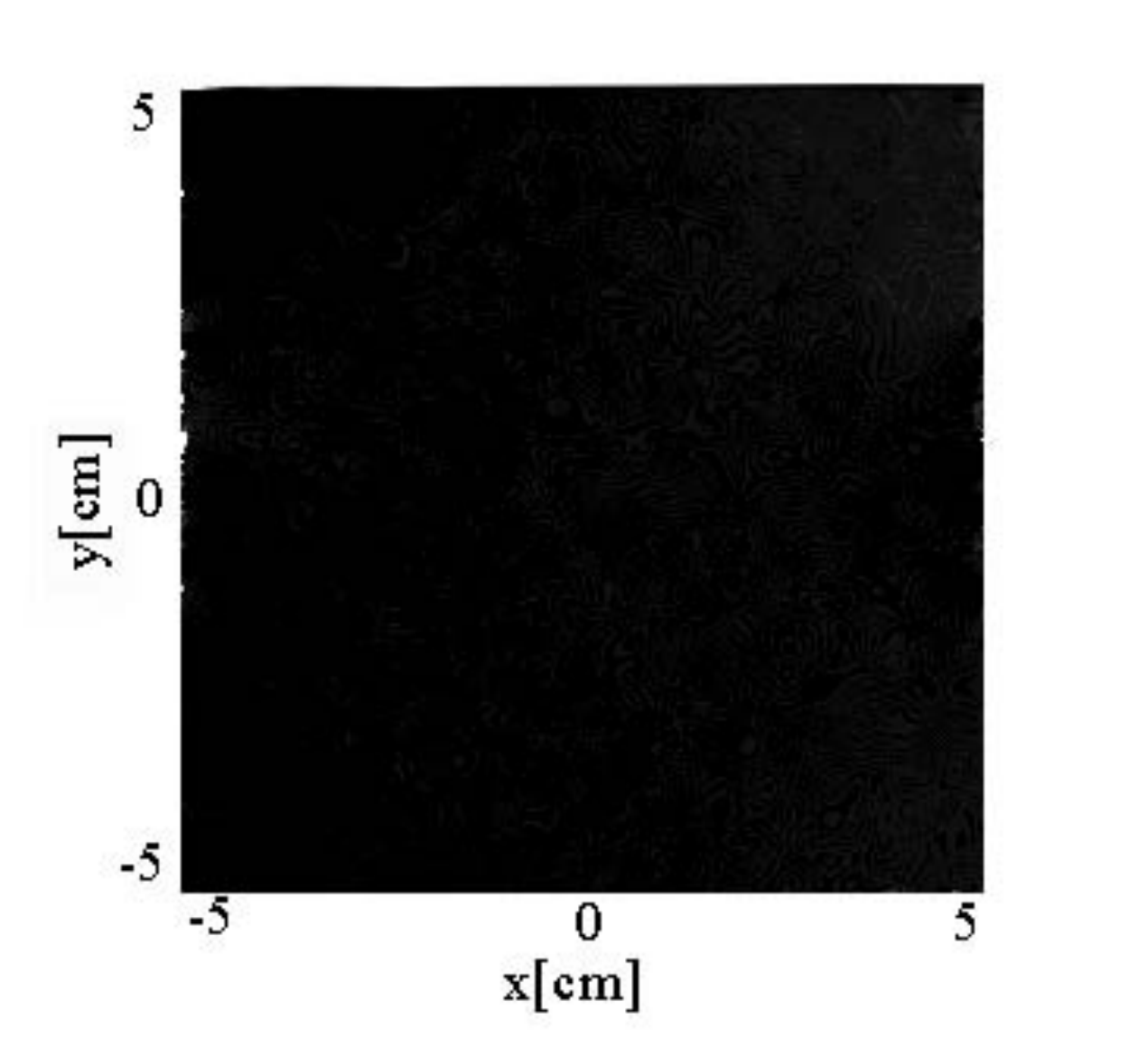}
\label{fig:tabsubfig31}
}
\caption{(Subplot (a)) Phase reconstructed obtained by means of CCS for SNR = 40 dB and $r=0.5$; (Subplot (b)) Phase reconstructed obtained by means of DCS for the same values of SNR and $r$.; (Subplot (c) and (d)) Corresponding error maps for CCS and DCS.}
\label{F5}
\end{figure}

Some typical phase reconstruction results are shown in Fig.~\ref{F5}, whose left and right subplots depict the phase estimates obtained using the CCS and DCS algorithms, respectively, for the case of $r = 0.5$. To visualize the differences more clearly, error maps for both methods are shown in subplot (c) and (d)(note for better visualization error map values are multiplied by factor of 10.). A close comparison with the original phase (as shown in subplot (a) of Fig.~\ref{F3}) reveals that DCS provides a more accurate recovery of the original $\phi$, which further supports the value of using the cross-derivative constraints. In fact, exploiting these constraints effectively amounts to using additional measurements, which are ignored in the case of CCS.

\begin{figure}[!t]
\centering
\includegraphics[width = 9cm, height = 6cm]{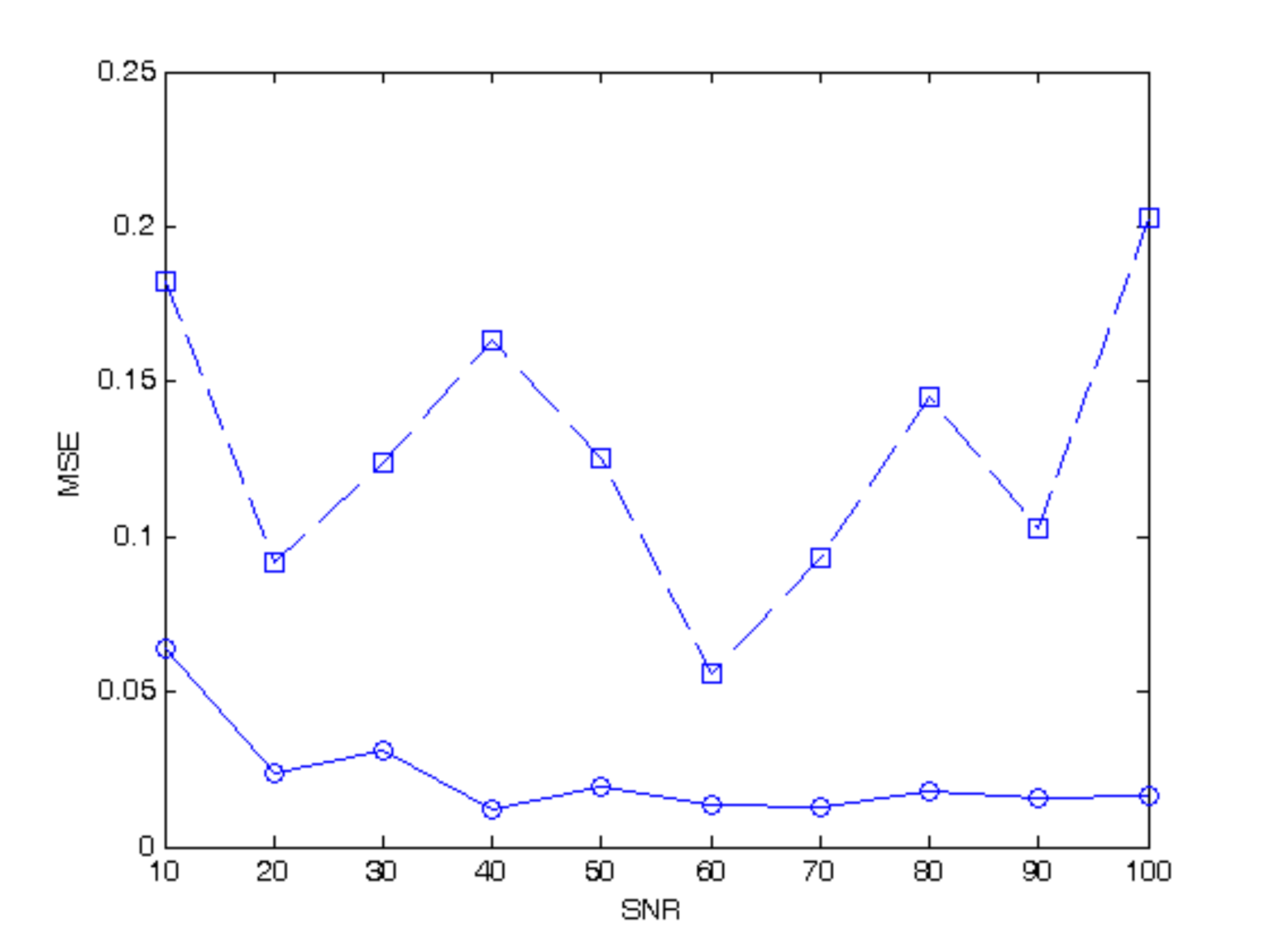}
\caption{MSE of phase reconstruction obtained with different methods as a function of SNR. Here, the dashed and solid lines correspond to CCS and DCS, respectively, and $r=0.5$.}
\label{F6}
\end{figure}

To investigate the robustness of the compared algorithms towards the influence of additive noises, their performances have been compared for a range of SNR values. The results of this comparison are summarized in Fig.~\ref{F6}. Since the cross-derivative constraints exploited by DCS effectively restrict the feasibility region for an optimal solution, the algorithm exhibits a substantially better robustness to the additive noise as compared to the case of CCS. This fact represents another beneficial outcome of incorporating the cross-derivative constraints in the process of phase recovery.

It should be taken into account that, although the shape of $\phi$ does not change the energy of the PSF $i$, it plays a crucial role in determining its spatial behaviour. In the section that follows, it will be shown that even small inaccuracies in reconstruction of $\phi$ could be translated into dramatic  difference in the quality of image deconvolution.

\subsection{Image deconvolution}
\begin{figure}[!t]
\subfigure[]{
\includegraphics[width = 4cm]{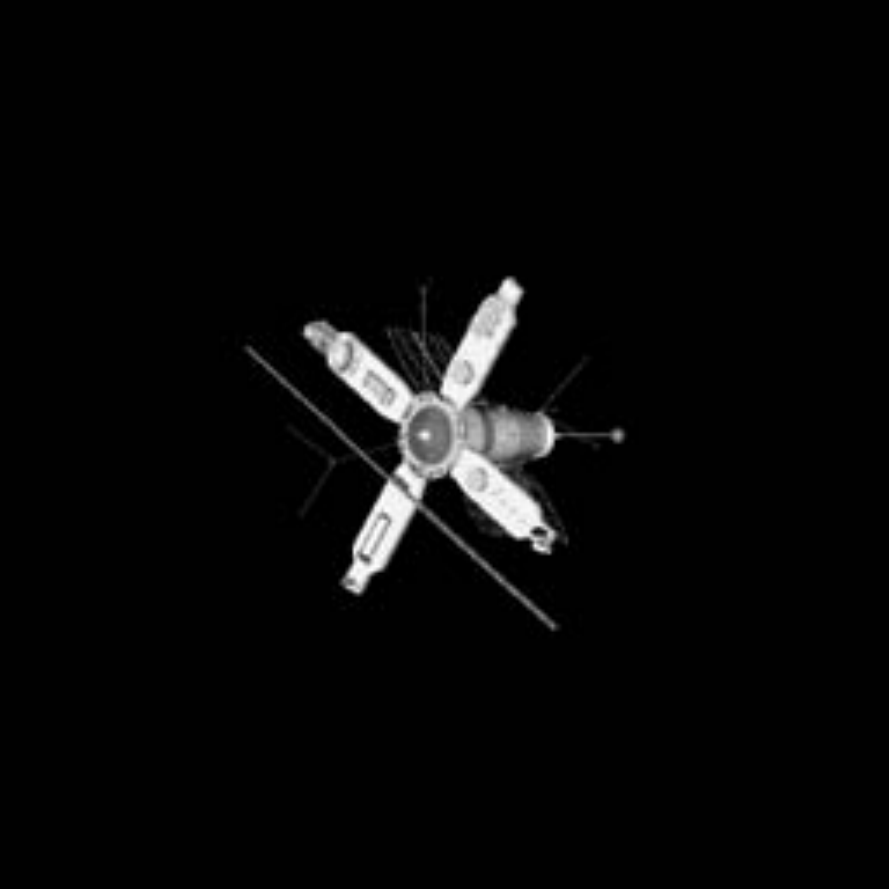}
\label{fig:tabsubfig2}
}
\subfigure[]{
\includegraphics[width = 4cm]{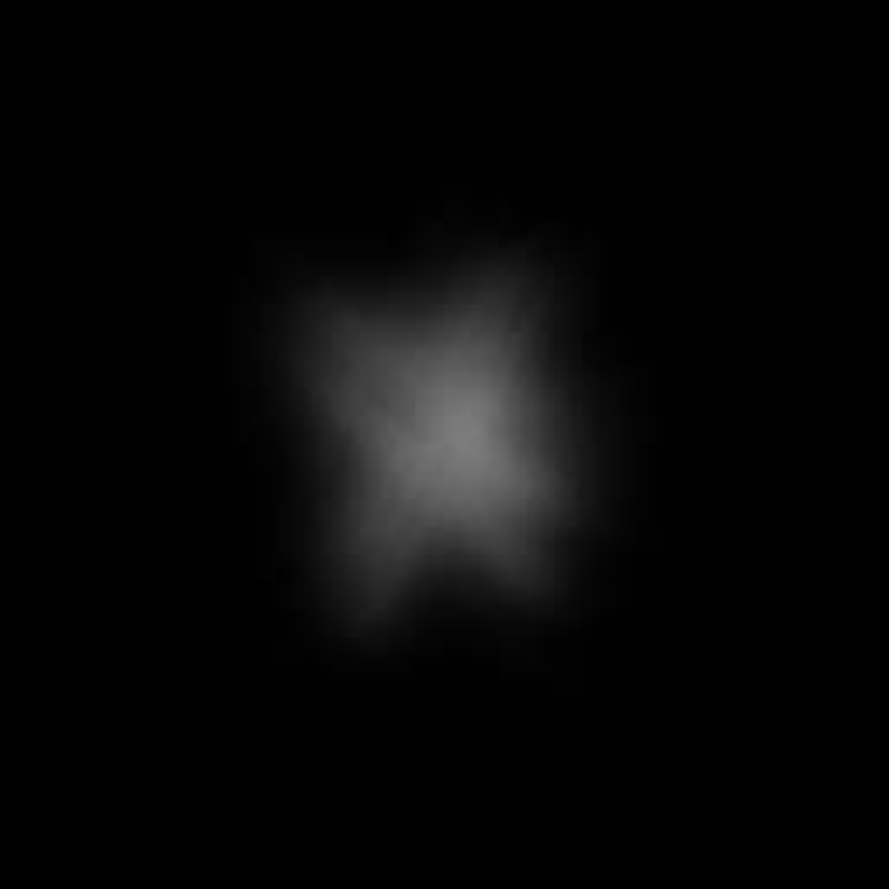}
\label{fig:tabsubfig3}
}
\caption{Satellite image (subplot (a)) and its blurred and noisy version (subplot (b)).}
\label{fig12}
\end{figure}

As a next step, the phase estimates obtained using the CCS- and DCS-based methods for $r=0.5$ were combined with the aperture function $A$ to result in their respective estimates of the PSF $i$. These estimates were subsequently used to deconvolve a number of test images such as ``Satellite", ``Saturn", ``Moon" and ``Galaxy". All the test images were blurred with an original PSF, followed by their contamination with additive Gaussian noise of different levels. As an example, Fig.~\ref{fig12} shows the ``Satellite" image (subplot (a)) along with its blurred and noisy version (subplot (b)).

\begin{figure}[t]
\subfigure[]{
\includegraphics[width = 4cm]{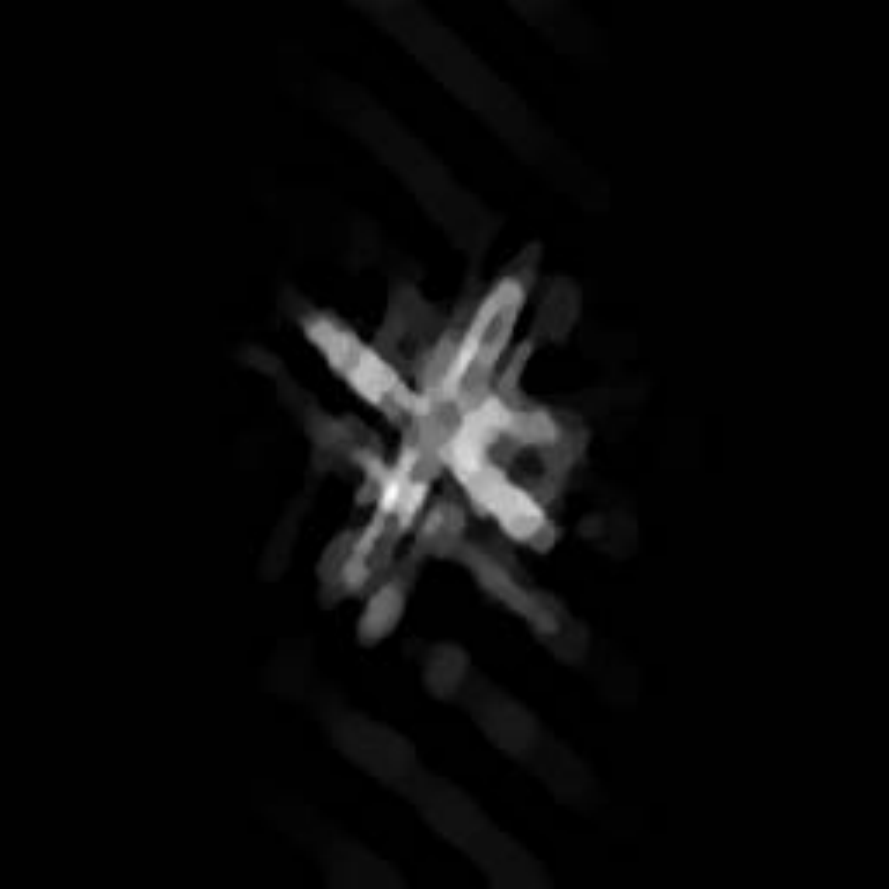}
\label{fig:tabsubfig2}
}
\subfigure[]{
\includegraphics[width = 4cm]{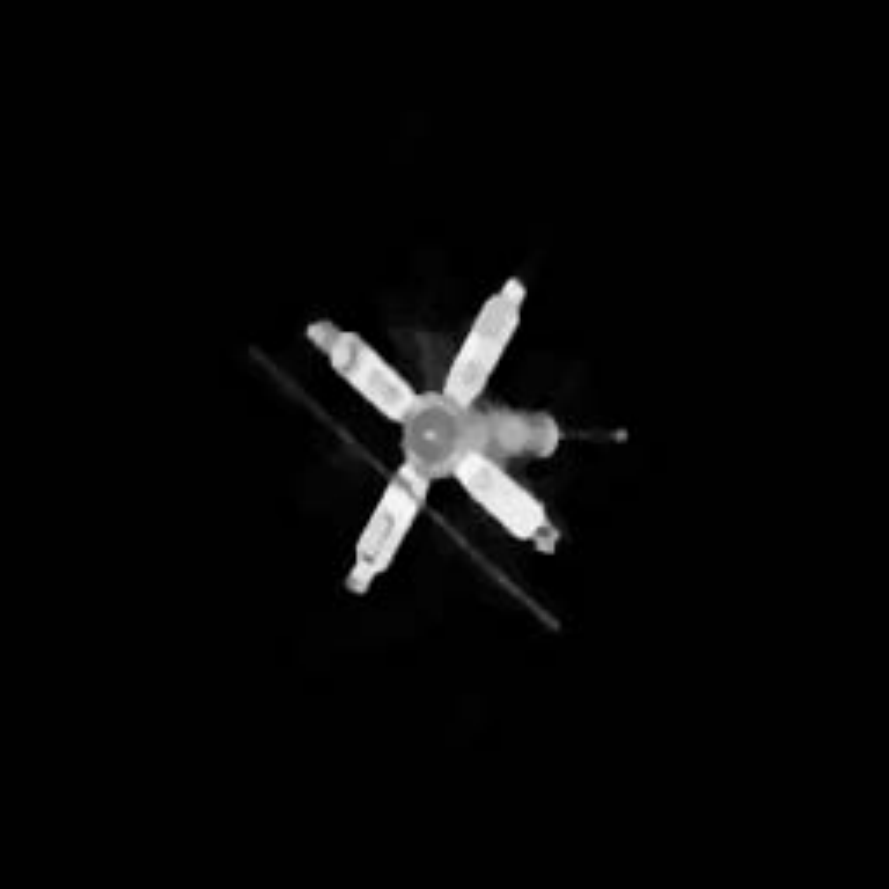}
\label{fig:tabsubfig3}
}
\caption{(Subplot (a)) Image estimate obtained with the CCS-based method for phase recovery (SSIM = 0.917); (Subplot (b)) Image estimate obtained with the DCS-based method for phase recovery (SSIM = 0.781).}
\label{fig11}
\end{figure}

Using the PSF estimates, the deconvolution was carried out using the method detailed in \cite{11}. For the sake of comparison, the deconvolution was also performed using the PSF recovered from dense sampling (DS) of $\phi$. Note that this reconstruction is expected to have the best accuracy, since it neither involves undersampling nor requires a CS-based phase estimation. All the deconvolved images have been compared with their original counterparts in terms of PSNR as well as of the structural similarity index (SSIM) of \cite{36}, , which is believed to be a better indicator of perceptual image quality \cite{100}. The resulting values of the comparison metrics are summarized in Table 1, while Fig.~\ref{fig11} shows the deconvolution results produced by the CCS- and DCS-based methods.

The above results demonstrate the importance of accurate phase recovery, where even a relatively small phase error can have a dramatic effect on the quality of image deconvolution. Under such conditions, the proposed method produces image reconstructions of a superior quality as compared to the case of CCS. Moreover, comparing the results of Table 1, one can see that DS only slightly outperforms DCS in terms of PSNR and SSIM, while in many practical cases, the difference between the performances of these methods are hard to detect visually.

\begin{table*}[t]
\centering \caption{SSIM and PSNR comparisons of phase recovery
results}
\label{tab:denresults}
\begin{tabular}{l|cccc|cccc|cccc|cccc}
\hline
Image & \multicolumn{4}{c}{Satellite} & \multicolumn{4}{c}{Saturn} & \multicolumn{4}{c}{Moon} & \multicolumn{4}{c}{Galaxy}\\
\hline
Noise std & $10^{-5}$ & $0.001$ & $0.003$ & $0.005$ & $10^{-5}$ &  $0.001$ & $0.003$ & $0.005$  & $10^{-5}$ &  $0.001$ & $0.003$ & $0.005$  & $10^{-5}$ &  $0.001$ & $0.003$ & $0.005$ \\
\hline
& \multicolumn{16}{c}{PSNR comparison (in dB)}\\
\hline
Blurred               & 14.06 & 14.06 & 14.06 & 14.05 & 17.78 & 17.78 & 17.78 &  17.78 & 19.98 & 19.97 & 19.97 &  19.97 & 18.79 & 18.79 & 18.78 &  18.78\\
DS                & 27.97 & 27.75 & 25.97 & 22.43 & 31.49 & 31.08 & 28.50 &  23.89 & 25.06 & 25.04 & 24.83 &  21.76 & 23.58 & 23.60 & 23.38 &  20.93\\
CS                    & 17.06 & 16.93 & 16.54 & 15.63 & 23.42 & 23.38 & 22.80 &  20.55 & 22.36 & 22.38 & 22.30 &  19.73 & 21.16 & 21.12 & 20.64 & 18.46\\
DCS                   & 27.42 & 27.22 & 25.56 & 22.22 & 31.02 & 30.65 & 28.30 &  23.72 & 25.00 & 24.99 & 24.78 &  21.73 & 23.52 & 23.54 & 23.32 & 20.86\\
\hline
& \multicolumn{16}{c}{SSIM comparison}\\
\hline
Blurred               & 0.200 & 0.200 & 0.199 & 0.197 & 0.226 & 0.226 & 0.226 & 0.175 & 0.512 & 0.512 & 0.509 & 0.504 & 0.257 & 0.257 & 0.257 & 0.254\\
DS               & 0.730 & 0.720 & 0.554 & 0.269 & 0.688 & 0.660 & 0.506 & 0.228 & 0.645 & 0.642 & 0.607 & 0.552 & 0.493 & 0.495 & 0.501 & 0.397\\
CS                    & 0.349 & 0.344 & 0.306 & 0.206 & 0.424 & 0.416 & 0.348 & 0.212 & 0.539 & 0.538 & 0.493 & 0.488 & 0.348 & 0.347 & 0.326 & 0.224\\
DCS                   & 0.674 & 0.667 & 0.519 & 0.263 & 0.656 & 0.641 & 0.483 & 0.223 & 0.643 & 0.640 & 0.604 & 0.549 & 0.490 & 0.491 & 0.501 & 0.393\\
\hline
\end{tabular}
\end{table*}

\section{Discussion and conclusions}
In the present paper, the applicability of DCS to the problem of reconstruction of optical images has been demonstrated. It was shown that, in the presence of atmospheric turbulence, the phase $\phi$ of the GPF $P = A \, e^{\jmath \phi}$ is a random function, which needs to be measured using adaptive optics. To simplify the complexity of the latter, a CS-based approach has been proposed. As opposed to CCS, however, the proposed method performs phase reconstruction subject to an additional constraint, which stems from the property of $\nabla \phi$ to be a potential field. The resulting algorithm (referred to as the DCS method) has been shown to yield phase estimates of substantially better quality as compared to the case of CCS.

In this paper, our main focus has been on simplifying the structure of the SHI through reducing the number of its wavefront lenses, while compensating for the effect of undersampling by using the theory of CS augmented by the cross-derivative constraint. The solution was computed using the Bregman algorithm, which provides a computationally efficient framework to carry out the constrained phase recovery. Moreover, the resulting phase estimates were used to recover their associated PSF, which was subsequently used for image deconvolution. It was shown that the DCS-based estimation of $\phi$ with $r=0.3$ results in image reconstructions of the quality comparable to that of DS, while substantially outperforming the results obtained with CCS.

While the proposed method offers a practical solution to the problem of phase estimation in adaptive optics, some interesting questions about the theoretical aspects of DCS still lay open. In particular, the question of theoretical performance of CS in the presence of side information on the source signal needs to be addressed through future research.

\section*{Acknowledgment}
This work was supported in part by the Natural Sciences and Engineering Research Council of Canada and in part by Ontario Early Researcher Award program, which are gratefully acknowledged. The authors would also like to acknowledge Sudipto Dolui for his helpful comments as well as for providing deconvolution codes.

\bibliography{DCS}
\bibliographystyle{IEEEbib}
\end{document}